\documentclass[12pt]{article}

\usepackage{amsmath,amsthm,amssymb,mathrsfs}
\usepackage[dvips]{graphicx}
\usepackage[dvips]{epsfig}
\usepackage{color}

\def\titlerunning#1{\gdef\titrun{#1}}
\makeatletter
\def\author#1{\gdef\autrun{\def\and{\unskip, }#1}\gdef\@author{#1}}
\def\address#1{{\def\and{\\\hspace*{18pt}}\renewcommand{\thefootnote}{}%
\footnote {#1}}%
\markboth{\autrun}{\titrun}}
\makeatother
\def\email#1{e-mail: #1}

\def\i{{\bf i}}
\def\j{{\bf j}}
\def\k{{\bf k}}

\def\p{{\bf p}}
\def\q{{\bf q}}
\def\r{{\bf r}}
\def\s{{\bf s}}

\def\u{{\bf u}}
\def\v{{\bf v}}
\def\x{{\bf x}}

\def\y{{\bf y}}

\def\Q{{\bf Q}}
\def\R{{\bf R}}
\def\S{{\bf S}}

\def\RR{{\mathop{{\rm I}\kern-.2em{\rm R}}\nolimits}}

\def\scal{{\rm {scal}}\,\!}
\def\vect{{\rm {vect}}\,\!}

\newcommand{\X}{{\cal X}}
\newcommand{\Y}{{\cal Y}}
\newcommand{\Z}{{\cal Z}}
\renewcommand{\P}{{\cal P}}
\renewcommand{\Q}{{\cal Q}}
\renewcommand{\R}{{\cal R}}
\renewcommand{\S}{{\cal S}}
\newcommand{\U}{{\cal U}}
\newcommand{\V}{{\cal V}}
\newcommand{\A}{{\cal A}}
\newcommand{\B}{{\cal B}}

\newcommand{\nn}{\mathbb{N}}
\newcommand{\rr}{\mathbb{R}}
\newcommand{\cc}{\mathbb{C}}
\newcommand{\hh}{\mathbb{H}}
\newcommand{\bbS}{\mathbb{S}}

\def\half{\textstyle\frac{1}{2}\displaystyle}

\newcommand{\be}{\begin{equation}}
\newcommand{\ee}{\end{equation}}
\newcommand{\ba}{\begin{eqnarray}}
\newcommand{\ea}{\end{eqnarray}}
\newcommand{\bi}{\begin{itemize}}
\newcommand{\ei}{\end{itemize}}

\newtheorem{dfn}{Definition}
\newtheorem{rmk}{Remark}
\newtheorem{lma}{Lemma}

\newtheorem{thm}{Theorem}

\newtheorem{exm}{Example}

\begin{document}

\titlerunning{}
\title{\bf Solution of a quadratic quaternion 
equation with mixed coefficients}

\author{Rida~T.~Farouki,
Graziano Gentili, Carlotta Giannelli, \\ 
Alessandra Sestini, and Caterina Stoppato}

\date{}

\maketitle

\address{Rida~T.~Farouki,
Department of Mechanical and Aerospace Engineering,
University of California, Davis, CA 95616--5294, USA,
\email{farouki@ucdavis.edu}
\and
Graziano Gentili, Alessandra Sestini,
Dipartimento di Matematica e Informatica \\ ``U. Dini,'' 
Universit\`a di Firenze, Viale Morgagni 67/A, I--50134 Firenze, Italy,
\email{gentili@math.unifi.it, alessandra.sestini@unifi.it}
\and
Carlotta Giannelli, Caterina Stoppato,
Istituto Nazionale di Alta Matematica, 
Unit\`a di Ricerca di Firenze c/o DiMaI ``U. Dini,'' 
Universit\`a di Firenze, Viale Morgagni 67/A, I-50134 Firenze, Italy,
\email{carlotta.giannelli@unifi.it, stoppato@math.unifi.it}
}

\begin{abstract}
A comprehensive analysis of the morphology of the solution space for a 
special type of quadratic quaternion equation is presented. This equation,
which arises in a surface construction problem, incorporates linear terms 
in a quaternion variable and its conjugate with right and left quaternion 
coefficients, while the quadratic term has a quaternion coefficient placed
between the variable and its conjugate. It is proved that, for generic 
coefficients, the equation has two, one, or no solutions, but in certain 
special instances the solution set may comprise a circle or a 3--sphere 
in the quaternion space $\hh$. The analysis yields solutions 
for each case, and intuitive interpretations of them in terms of the 
four--dimensional geometry of the quaternion space $\hh$.
\end{abstract}

\newpage

%%%%%%%%%%%%%%%%%%%%%%%%%%%%%%%%%%%%%%%%%%%%%%%%%%%%%%%%%%%%%%%%%%%%%%%%%%%%

\section{Introduction}

The real quaternions, discovered by Hamilton in 1843, form the first known 
algebra that involves a non--commutative product, denoted by $\hh$. This 
property makes the solution of equations featuring quaternion coefficients 
and unknowns a much subtler and richer problem than in the case of real or 
complex numbers. The present study is devoted to developing a comprehensive 
solution for a novel type of quadratic equation, of the form
\be
\label{Xeqn}
\X\,\P\X^* \,+\, \X\Q\,+\, \R\,\X^* \,=\, \S \,,
\ee
in which both the coefficients $\P,\Q,\R,\S$ and the variable $\X$ are 
quaternions ($\X^*$ being the conjugate of $\X$). Since we are only 
interested in the genuinely quadratic case of equation~\eqref{Xeqn}, 
we assume henceforth that $\P \neq 0$. We further emphasize that only the 
algebra $\hh $ of \emph{real} quaternions is considered herein: the results 
do not apply, for example, to the complexified algebra $\hh\otimes \cc$. 

A complete analysis of equation~\eqref{Xeqn} reveals that, in addition 
to cases with (at most) two distinct \emph{point} solutions, special 
values of the coefficients yield singly--infinite (\emph{circular}) and 
triply--infinite (\emph{3--sphere}) families of solutions. A proper
identification and treatment of these degenerate cases is therefore an 
essential feature of any comprehensive solution procedure.

The motivation for studying equation \eqref{Xeqn} arises \cite{farouki15} 
from the construction of a surface patch $\x(u,v)$ defined on $(u,v)\in
[\,0,1\,]^2$ with prescribed boundary curves, such that the $v=\mbox{constant}$
 isoparametric curves are all polynomial \emph{Pythagorean--hodograph} (PH) 
\emph{curves} \cite{farouki08}. A brief review of the surface construction 
problem can be found in Section~\ref{sec:phcurves} below. At present, we 
highlight some unusual properties of equation \eqref{Xeqn} that distinguish 
it from prior studies of quadratic (and higher--order) quaternion equations.

The earliest investigations of quaternion equations are found in the papers 
of Niven \cite{niven41,niven42}, Eilenberg and Niven \cite{eilenberg44}, and 
Gordon and Motzkin \cite{gordon65}. A special case is that of equations of 
the form $f(\X)=0$, where $f(\X)$ is an element of the algebra of quaternion 
polynomials $\hh[\X]$. In this algebra, a convention is fixed for the relative 
position of the quaternion coefficients and powers of the quaternion variable 
in each monomial. Thus, one speaks of a left/right quaternion polynomial if 
the coefficients all appear to the left/right of the corresponding powers 
of the quaternion unknown. 

Several studies \cite{farouki13,huang02,jia09} have been specifically
concerned with (monic) quadratic quaternion equations specified by unilateral 
coefficients. In this case an essentially closed--form solution, requiring a
determination of the positive real root of a real cubic equation by Cardano's 
method, is possible (including a complete enumeration of special cases).
Although no closed--form solution for the roots of higher--order polynomials
is possible, there has been considerable progress in elucidating their 
fundamental nature, and in developing numerical methods to compute them
\cite{damiano10,deleo06,gentili08a,gentili08b,gentili08c,janovska10,kalantari13,pogorui04,pumplun02,serodio01a,serodio01b,topuridze09}.
In particular, it has been shown that the set of roots of any polynomial 
$f(\X) \in \hh[\X]$ is a finite union of singletons and $2$--spheres.

On the other hand, one may consider the class $\mathbf{E}$ of functions 
$f:\hh\to\hh$ defined as finite sums of monomials of type $\A_0\X\A_1\cdots
\X\A_n$ with $\A_0,\ldots,\A_n \in \hh$ and $n \in \nn$. This class is 
extremely large. In fact, it includes each of the functions that map any 
quaternion $\X=x_0+x_1\i+x_2\j+x_3\k$ to one of its four real components 
$x_0,\ldots,x_3$ since, by direct computation, 
\begin{align*}
x_0 &= \frac 14 (\X -\i\,\X\,\i-\j\,\X\,\j-\k\,\X\,\k) , \\
x_1 &= \frac \i 4 (-\X +\i\,\X\,\i-\j\,\X\,\j-\k\,\X\,\k) , \\
x_2 &= \frac \j 4 (-\X -\i\,\X\,\i+\j\,\X\,\j-\k\,\X\,\k) , \\
x_3 &= \frac \k 4 (-\X -\i\,\X\,\i-\j\,\X\,\j+\k\,\X\,\k) .
\end{align*}
Consequently, the class of equations $f(\X)=0$ with $f \in \mathbf{E}$ 
coincides with the class of all systems of four real polynomial equations 
in the four real variables $x_0,\ldots,x_3$. Moreover, $\mathbf{E}$ includes 
the map $\X \mapsto \X^*$, so equation \eqref{Xeqn} also has the form 
$f(\X)=0$ with $f \in \mathbf{E}$.

In the present paper, we show that the subclass of equations of the form 
\eqref{Xeqn} is special enough to allow a comprehensive solution. We show 
that generically there are two, one, or no solutions, which may be determined 
geometrically. The remaining special cases, in which an entire circle or 
$3$--sphere of solutions may occur, are also studied in detail. The method 
adopted herein is to reduce all instances of equation~\eqref{Xeqn} to two 
special cases, which may be treated as systems of real quadratic equations 
--- as one may expect from the preceding discussion concerning the class 
of functions $\mathbf{E}$.

The plan for this paper is as follows. First, some preparatory notations 
and results are presented in Section~\ref{sec:notation}, and the surface 
construction problem resulting in equation~\eqref{Xeqn} is briefly reviewed 
in Section~\ref{sec:phcurves}. A comprehensive solution procedure,
including treatment of all special cases, is then developed in Section~\ref
{sec:solution}, and examples are presented to illustrate the different 
solution morphologies that may arise, depending on the nature of the 
coefficients. Finally, Section~\ref{sec:close} assesses the principal 
results of the present study.

%%%%%%%%%%%%%%%%%%%%%%%%%%%%%%%%%%%%%%%%%%%%%%%%%%%%%%%%%%%%%%%%%%%%%%%%%%%%

\section{Notations and preliminary results}
\label{sec:notation}

Before proceeding to the study of equation~\eqref{Xeqn}, we fix some 
notations and briefly recall the definition and basic properties of the 
algebra of quaternions $\hh$. It is the vector space $\rr^4$ endowed with 
a multiplicative operation through the following steps:
\bi
\item denote the standard basis of $\rr^4$ as $1,\i,\j,\k$;
\item let the (left or right) multiplication by $1$ have no effect on 
$1,\i,\j,\k$;
\item set $\i^2=\j^2=\k^2=-1$ and $\i\,\j=-\j\,\i=\k$, 
$\j\,\k=-\k\,\j=\i$, $\k\,\i=-\i\,\k=\j$;
\item extend the multiplication to all quaternions $\X=x_01+x_1\i+x_2\j+x_3\k$ 
(with $x_0,\ldots,x_3 \in \rr$) in a bilinear fashion.
\ei
The resulting algebra is associative and non--commutative. It is also
unitary, with identity element $1$. In view of this, we henceforth write 
$x_0$ instead of $x_0\,1$, and we identify the axis $\rr\,1$ with the real 
field $\rr$.

Moreover, $\hh$ is a skew field --- namely, every non--zero element admits 
a multiplicative inverse. In order to give a formula for the inverse, we 
introduce some further notations. For each quaternion $\X=x_0+x_1\i+x_2\j
+x_3\k \in \hh$, we denote by
\[
\x = x_1\i+x_2\j+x_3\k
\]
its \emph{vector} \emph{part}, and by $x_0$ its \emph{scalar} \emph{part}. 
The notations $\scal(\X)=x_0$ and $\vect(\X)= x_1\i+x_2\j+x_3\k$ are also 
used. A quaternion whose scalar part vanishes is called a \emph{pure vector} 
quaternion. To each quaternion $\X=x_0+\x=x_0+x_1\i+x_2\j+x_3\k$, one 
associates the \emph{conjugate}
\[
\X^* = x_0 - \x =x_0-x_1\i-x_2\j-x_3\k,
\]
such that $\X\,\X^*=x_0^2+|\x|^2 = x_0^2+x_1^2+x_2^2+x_3^2$ coincides with 
the square of the Euclidean norm $|\X|$ of $\X$, i.e.,
\[
\X\,\X^*=\X^*\X= |\X|^2.
\]
By analogy with the complex case, the (non--negative) real number $|\X|$ 
is also called the \emph{modulus} of $\X$. Clearly $|\X|^2$ is a positive 
real number for all non--zero $\X \in \hh$, and the preceding equation 
implies that $\X^*|\X|^{-2} =|\X|^{-2} \X^*$ is the (left and right) inverse 
of $\X$, i.e.,
\[
\X^{-1} = \X^*|\X|^{-2} =|\X|^{-2} \X^*.
\]
The following properties will also prove useful in our computations.
\bi
\item If we denote by $\langle\;,\,\rangle$ and $\times$, respectively, 
the Euclidean scalar product and the vector product in $\rr^3$, which we 
may identify with $\vect(\hh) = \rr \i + \rr \j + \rr \k$, then for all 
$\x,\y \in \vect(\hh)$, we have 
\[
\x\,\y = -\, \langle\x,\y\rangle + \x\times\y .
\]
\item If we denote by $\langle\;,\,\rangle$ the Euclidean scalar product 
in $\rr^4 = \hh$, then for all $\X,\Y \in \hh$ we have
\[
\langle\X,\Y\rangle 
= \scal (\X\,\Y^*) .
\]
\ei
Finally, we denote by $\bbS = \{\X : \X^2 = -1\} = \{\x \in \vect(\hh) : 
|\x| =1\}$ the $2$--sphere of \emph{pure vector quaternion units}. Then the
following properties hold.
\bi
\item For any fixed $\u \in \bbS$ the subalgebra of $\hh$ generated by 
$1$ and $\u$, namely $\rr+\u\,\rr$, is isomorphic to the complex plane and 
is denoted by $\cc_\u$.
\item $\cc_\u=\cc_\v$ if and only if $\u = \pm \v$; otherwise $\cc_\u\cap
\cc_\v = \rr$.
\item Every $\X \in \hh\setminus \rr$ (i.e., $\X=x_0+\x$ with $\x\neq{\bf 0}$) 
belongs to $\cc_{\hat\x}$, where 
\[
\hat\x:= \frac{\x}{|\x|}.
\]
Every $\X \in \rr$ belongs to $\cc_\u$ for all $\u \in \bbS$. Therefore:
\[
\hh = \bigcup_{\u \in \bbS}\cc_\u.
\]
\ei
These properties will be used extensively in the proofs of our main results.

%%%%%%%%%%%%%%%%%%%%%%%%%%%%%%%%%%%%%%%%%%%%%%%%%%%%%%%%%%%%%%%%%%%%%%%%%%%%

\section{Surface construction problem}
\label{sec:phcurves}

The motivation for the study of equation~\eqref{Xeqn} arises \cite
{farouki15} in the construction of a surface patch $\x(u,v)$ for 
$(u,v)\in[\,0,1\,]^2$ with prescribed boundary curves, such that the 
$v=\mbox{constant}$ isoparametric curves are Pythagorean--hodograph 
(PH) curves.\footnote{The PH curves have rational unit tangents, 
polynomial arc length functions, and many other attractive properties 
--- complete details may be found in \cite{farouki08}.} Such a surface 
is obtained by integrating the expression
\be 
\label{xu}
\x_u(u,v) \,=\, \A(u,v)\,\i\,\A^*(u,v) \,,
\ee
where $\A(u,v)$ is a bivariate tensor--product quaternion polynomial
\be 
\label{A}
\A(u,v) \,=\, \sum_{i=0}^m \sum_{j=0}^n 
\A_{ij}\,b^m_i(u)b^n_j(v) \,. 
\ee 
expressed in terms of the Bernstein basis 
\[
b^d_k(t) \,=\, \binom{d}{k}(1-t)^{d-k}t^k \,, \quad k=0,\ldots,d \,.
\]
The simplest non--trivial solutions correspond to $m=n=2$. Integrating 
\eqref{xu} allows the surface to be expressed in terms of \emph{B\'ezier 
control points} $\p_{ij}$ \cite{farin} as
\be 
\label{x} 
\x(u,v) \,=\, \sum_{i=0}^5 \sum_{j=0}^4
\p_{ij}\,b^5_i(u)\,b^4_j(v) \,.
\ee 
The points $\p_{ij}$ for $i>0$ can be expressed in terms of the coefficients 
${\cal A}_{ij}$ for $0 \le i,j \le 2$, while $\p_{0j}$ for $j=0,\ldots,4$ 
amount to free integration constants that specify $\x(0,v)$. Note, in 
particular, that $\p_{i0}$ and $\p_{i4}$ for $i=0,\ldots,5$ depend only on 
the coefficients ${\cal A}_{i0}$ and ${\cal A}_{i2}$ for $i=0,\ldots,2$, 
respectively.

Now $\A_{00},\A_{10},\A_{20}$ and $\A_{02},\A_{12},\A_{22}$ may be used to 
fix the boundary PH curves $\x(u,0)$ and $\x(u,1)$ --- i.e., to determine 
$\p_{i0}$ and $\p_{i4}$ for $i=0,\ldots,5$ --- as Hermite interpolants \cite
{farouki08b}. The remaining coefficients $\A_{01},\A_{11},\A_{21}$ must 
then be used to achieve desired positions for the three interior control 
points $\p_{51}$, $\p_{52}$, $\p_{53}$ that specify $\x(1,v)$. This problem is 
under--determined, with three free parameters, but we assume that they are 
chosen \emph{a priori}. One can formulate $\p_{51}$, $\p_{52}$, $\p_{53}$ as 
quadratic expressions in $\A_{01},\A_{11},\A_{21}$ and known quantities. Two 
of the resulting equations are linear in $\A_{01}$, $\A_{21}$ and can be 
used to express these unknowns in terms of $\A_{11}$. Finally, substituting 
these expressions into the third equation yields the form \eqref{Xeqn}, 
where $\P,\Q,\R,\S$ are known, and we set $\X=\A_{11}$. Complete details 
may be found in \cite{farouki15}.

%%%%%%%%%%%%%%%%%%%%%%%%%%%%%%%%%%%%%%%%%%%%%%%%%%%%%%%%%%%%%%%%%%%%%%%%%%%%

\section{The solution procedure}
\label{sec:solution}

Different 
approaches to equation~\eqref{Xeqn} are appropriate, according to whether 
or not the coefficient $\P$ lies on the real axis $\rr$. These cases will 
be analyzed in detail in Sections~\ref{sec:real} and \ref{sec:nonreal} 
below. The principal results of this analysis are sketched in the 
following theorem.

\begin{thm}\label{general}
For generic $\P,\Q,\R,\S \in \hh$ with $\P \neq 0$, equation~\eqref{Xeqn}
has two, one, or no solutions. However, there are two non--generic instances
of $\P,\Q,\R,\S$ in which equation~\eqref{Xeqn} admits infinitely--many 
solutions --- namely, a $3$--sphere of solutions 
or a circle of solutions.
\end{thm}

\medskip\noindent
We now proceed to investigate individually the two cases $\P \in \rr$ 
and $\P \not \in\rr$ of equation~\eqref{Xeqn}, resulting in Theorems~\ref
{thm:real} and \ref{thm:nonreal} below --- which immediately imply 
Theorem~\ref{general}.

\subsection{The case $\P \in \rr$}
\label{sec:real}

We begin with the special case of equation~\eqref{Xeqn} in which the 
coefficient $\P$ lies on the real axis $\rr$. Since we will have occasion 
to refer individually to the scalar and vector parts of the quaternion 
variable and coefficients, we recall from Section~\ref{sec:notation} their 
splitting denoted by $\X=x_0+\x$, $\P=p_0+\p$, $\Q=q_0+\q$, $\R=r_0+\r$, 
$\S=s_0+\s$. 

\begin{thm}
\label{thm:real}
If $\P=p_0 \in \rr$ with $p_0 \neq 0$, equation~\eqref{Xeqn}
becomes
\be
\label{Xeqn2}
p_0|\X|^2 \,+\, \X\Q\,+\, \R\,\X^* \,=\, \S,
\ee
whose solutions in $\hh$ are the points $\X=\Y-(\Q^*+\R)/2p_0$, where $\Y$
satisfies
\be
\label{sphere}
|\Y|^2=\rho, \quad 
\rho:=\frac{\scal(\S)}{p_0} +\frac{|\Q+\R^*|^2}{4 p_0^2}
\ee
and
\be
\label{line}
\vect(\Y\,(\Q-\R^*))\,=\, \vect\left(\S+ \frac{\R\Q}{p_0}\right).
\ee
In particular, we may distinguish the following cases.
\begin{enumerate}
\item When $\rho<0$, there is no solution. 
\item When $\rho=0$, either $\X=-(\Q^*+\R)/2p_0$ is the unique solution or 
there is no solution, depending on whether or not
\[
\S+ \frac{\R\Q}{p_0}\in \rr.
\]
\item When $\rho>0$, then
\begin{enumerate}
\item if $\Q=\R^*$ the set of solutions is the $3$--sphere~\eqref{sphere} 
or the empty set, depending on whether or not $\S \in \rr$;
\item if $\Q\neq\R^*$ there are two, one, or no solutions, namely, the 
points $\X = \Y - (\Q^*+\R)/2p_0$ where $\Y$ is a point of intersection 
of the $3$--sphere~\eqref{sphere} with the affine line~\eqref{line}. 
Specifically,
\[
\Y=\left[\vect\!\left(\S+ \frac{\R\Q}{p_0}\right)\pm \sqrt{\Delta}\,\right] 
(\Q-\R^*)^{-1}
\]
where $\Delta$ is given by 
\[
\Delta = \rho\,|\Q-\R^*|^2-
\left|\vect\!\left(\S+ \frac{\R\Q}{p_0}\right)\right |^2.
\]
A positive, zero, or negative $\Delta$ identifies cases with two, one, 
or no solutions.
\end{enumerate}
\end{enumerate}
\end{thm}

\proof
Setting $\Y:=\X+(\Q^*+\R)/2p_0$ and substituting $\X=\Y-(\Q^*+\R)/2p_0$ 
into equation~\eqref{Xeqn2}, we obtain
\begin{align*}
\S =\ & p_0\left|\Y - \frac{\Q^*+\R}{2p_0}\right|^2 + 
\left(\Y - \frac{\Q^*+\R}{2p_0}\right)\Q + 
\R \left(\Y - \frac{\Q^*+\R}{2p_0}\right)^* \\
=\ & p_0|\Y|^2- \Y\frac{\Q+\R^*}{2} + \Y\Q 
- \frac{\Q^*+\R}{2} \Y^* + \R \Y^* \\
& + \frac{|\Q^*+\R|^2}{4p_0} - \frac{\Q^*+\R}{2p_0} \Q 
- \R \frac{\Q+\R^*}{2p_0} \\
=\ & p_0|\Y|^2+ \Y\frac{\Q-\R^*}{2} +\frac{-\Q^*+\R}{2} \Y^* \\
&+ \frac{|\Q^*+\R|^2}{4p_0} - \frac{|\Q^2|+2\R\Q+|\R|^2}{2p_0}\\
=\ & p_0|\Y|^2+ \Y \tilde \Q -\tilde \Q^* \Y^*
- \frac{|\Q^*+\R|^2}{4p_0} +\frac{-\R\Q+\Q^*\R^*}{2p_0}
\end{align*}
where $\tilde\Q=\half(\Q-\R^*)$. Hence, setting
\begin{align*}
\tilde \S =\ & \S +\frac{|\Q^*+\R|^2}{4p_0} + \frac{\R\Q-\Q^*\R^*}{2p_0} ,
\end{align*}
we have
\be
\label{realreduction}
p_0|\Y|^2+ \Y \tilde \Q -\tilde \Q^* \Y^* = \tilde \S .
\ee
Writing $\tilde \S = \tilde s_0 + \tilde \s$, equation~\eqref{realreduction} 
is equivalent to the system
\be
\label{realsystem}
p_0\, |\Y|^2 = \tilde s_0 , \qquad
\Y \tilde \Q -\tilde \Q^* \Y^* \,=\, \tilde \s .
\ee
The solutions to this system are those points that simultaneously satisfy 
$|\Y|^2=\rho$, with
\[ 
\rho = \frac{\tilde s_0}{p_0} = \frac{s_0}{p_0}+\frac{|\Q^*+\R|^2}{4p_0^2},
\] 
and 
\be
\label{affinespace}
\vect(\Y(\Q-\R^*)) \,= \, 2\, \vect(\Y \tilde \Q) \,=\, \tilde \s\, =\, \vect\left(\S + \frac{\R\Q}{p_0}\right).
\ee
 If $\rho<0$, then there is 
no solution in $\hh$, whence case {\it 1}.  If $\rho=0$, then $0$ is the 
only solution to $|\Y|^2=\rho$. It is also a solution of equation~\eqref
{affinespace} if, and only if, $\tilde\s={\bf 0}$. This verifies case 
{\it 2}. 

Consider now case {\it 3}, in which $\rho > 0$. If $\Q = \R^*$,  the 
solution to~\eqref{affinespace} is either the 
entire space $\hh$ or the empty set, according to whether or not $\s=0$.
On the other hand, the solutions to 
\eqref{affinespace} comprise an affine line if  $\Q 
\neq \R^*$. In this case, the solutions to \eqref{realsystem} are the points
\[
\Y = \left( \xi+\S+\frac{\R\Q}{p_0} \right)(\Q-\R^*)^{-1} ,
\]
where $\xi\in\rr$ is a root of the real quadratic equation
\[
\left|\, \xi+\S+\frac{\R\Q}{p_0} \,\right|^2 = \rho\,|\Q-\R^*|^2 .
\]
A positive, zero, or negative discriminant leads to two, one, 
or no solutions.

Finally, the translation $\X= \Y- (\Q^*+\R)/2p_0$ yields the solutions 
of equation~\eqref{Xeqn} in terms of those of \eqref{realsystem}.
\qed

\bigskip

We now describe some examples that serve to illustrate all the possible 
types of solution sets to equation \eqref{Xeqn} covered by Theorem~\ref
{thm:real}.
We begin with two examples that illustrate cases {\it 1} and {\it 2} of
Theorem~\ref{thm:real}.

\begin{exm}
{\rm If $\P=1=-\S$ and $\Q=\R=0$, equation~\eqref{Xeqn} becomes
\[
|\X|^2+1=0,
\]
which clearly has no solution in $\hh$.}
\end{exm}

\begin{exm}
{\rm If $\P=\Q=\R=1=-\S$, equation~\eqref{Xeqn} becomes
\[
|\X|^2+\X+\X^*+1=0,
\]
which is equivalent to $|\X+1|^2=0$. Thus, $\X=-1$ is the only solution 
in $\hh$.}
\end{exm}

\medskip\noindent
The following two examples are instances of case {\it 3(a)} in Theorem~\ref
{thm:real}.

\begin{exm}
{\rm If $\P=\Q=\R=1$ and $\S=0$, equation~\eqref{Xeqn} becomes
\[
|\X|^2+\X+\X^*=0,
\]
which reduces to $|\X+1|^2=1$. The set of solutions in $\hh$ is therefore the 
$3$--sphere of radius $1$ centered at $-1$.}
\end{exm}

\begin{exm}
{\rm If $\P=\Q=\R=1$ and $\S = \i$, equation~\eqref{Xeqn} becomes
\[
|\X|^2+\X+\X^* = \i .
\]
Since this cannot be satisfied by any $\X$, there are no solutions in $\hh$.}
\end{exm}

\medskip\noindent
Next is a family of examples corresponding to case {\it 3(b)} in Theorem~\ref
{thm:real}.

\begin{exm}
{\rm If $\P=\Q=1=-\R$ and $\S=1+\s$, equation~\eqref{Xeqn} becomes
\[
|\X|^2+\X-\X^*=1+\s.
\]
The solutions correspond to the intersections of the $3$--sphere $|\X|^2=1$ 
with the line $\X-\X^*=\s$. If $\s=2\sin\theta\,\hat\s$, there are two or 
one intersections, namely $\pm \cos(\theta)+\sin(\theta)\,\hat\s$. If $|\s|
>2$, on the other hand, there is no intersection.}
\end{exm}

\subsection{The case $\P \not \in \rr$}
\label{sec:nonreal}

Consider now equation \eqref{Xeqn} when $\P$ is not a real number. In 
this case, we shall make use of a known result (see Section~2.1 of \cite
{farouki08b}) concerning the instance of \eqref{Xeqn} in which $\P$, $\S$ 
are pure imaginary quaternions $\p$, $\s$ and $\Q=\R=0$. Recall from 
Section~\ref{sec:notation} that $\hat\p:=\p/|\p|$ and $\hat\s:=\s/|\s|$ 
are unit vectors in the direction of $\p$ and $\s$, and $\cc_{\hat\p}$ is 
the 2--plane spanned by 1 and $\hat\p$.

\begin{dfn}
\label{plane}
Let $\p,\s$ be non--zero pure vector quaternions. If $\hat\s\neq-\hat\p$, 
the 2--plane $\Pi_{\p,\s}$ is defined by
\[
\Pi_{\p,\s} := \u_{\p,\s}\,\cc_{\hat\p} \,, \quad {\rm with} \;\; 
\u_{\p,\s} := \frac{\hat\p+\hat\s}{\left|\hat\p+\hat\s\right|} .
\]
On the other hand, if $\hat\s=-\hat\p$, then $\Pi_{\p,\s}$ is defined to 
be the $2$--plane through the origin that is orthogonal to $\hat\p$ in the 
$3$--space $\vect(\hh)$.
\end{dfn}

\medskip\noindent
Note that the $2$--plane $\Pi_{\p,\s}$ bisects the angle between $\hat\p$ 
and $\hat\s$. 
\begin{lma}
\label{circle}
Let $\p,\s$ be pure vector quaternions with $\p \neq {\bf 0}$. Then if
$\s\ne{\bf 0}$, the 2--plane $\Pi_{\p,\s}$ is the space of solutions $\V$ 
to equation 
\be
\label{bisector}
\V\,\hat\p-\hat\s\, \V = 0.
\ee 
Moreover, the set of solutions to the quaternion equation
\be
\label{imaginaryequation}
\X\,\p\,\X^* \,=\, \s
\ee
is the circle $\mathscr{C}_{\p,\s}$ in the 2--plane $\Pi_{\p,\s}$ with center 
$0$ and radius $ \tau = \sqrt{|\s|/|\p|}$. On the other hand, when $\s=0$, it 
is simply $\mathscr{C}_{\p,0}:=\{0\}$.
\end{lma}

\medskip\noindent
Although this result is already known \cite{farouki08b}, a proof is included 
below to make the presentation self--contained.

\proof
Clearly, if $\s=0$, the unique solution is $\X=0$. We therefore consider
only the case $\s\neq 0$. Equation~\eqref{imaginaryequation} then implies 
that
\[
|\X|^2 = \frac{|\s|}{|\p|}.
\] 
Therefore, any solution must have the form $\X=  \tau\,\U$ with $\tau^2=
|\s|/|\p|$ and $|\U|=1$. Substituting into \eqref{imaginaryequation} then 
gives $\U\,\hat\p\,\U^*=\hat\s$, whose solution set is the intersection 
of the $3$--sphere $|\U| = 1$ with the set of solutions $\V$ to \eqref
{bisector}. To verify Lemma~\ref{circle}, we show that this set of 
solutions is precisely $\Pi_{\p,\s}$.

If $\hat\s\neq-\hat\p$ it is evident by inspection that $\V=\u_{\p,\s}$ is 
a solution of \eqref{bisector}. Moreover, if $\V$ satisfies \eqref{bisector}, 
the product of $\V$ with any element of $\cc_{\hat\p}$ also satisfies it. 
Thus, $\Pi_{\p,\s}=\u_{\p,\s}\,\cc_{\hat\p}$ belongs to the space of solutions 
of \eqref{bisector}. We will now prove that the dimension of this space 
cannot exceed two. Indeed, an orthonormal transformation gives 
$\hat\p=\i$ and $\hat\s=\i\cos\varphi+\j \sin\varphi$  for some 
$\varphi \in \rr$ and the rank of the linear map 
$\V \mapsto \V\,\i-(\i\cos\varphi+\j\sin\varphi)\,\V$ cannot be less 
than $2$ since, by inspection, its values at $\V=\i$ and $\V=\k$ are linearly 
independent. Hence, the solution space of \eqref{bisector} coincides 
precisely with $\Pi_{\p,\s}$.

Otherwise, if $\hat\s=-\hat\p$, the solutions $\V=v_0+\v$ must satisfy
\[
0 = \V\,\hat\p+\hat\p\, \V = 2v_0\,\hat\p + 
\v\,\hat\p+\hat\p\,\v = 2v_0\,\hat\p - 2\,\langle\v,\hat\p\rangle
\]
i.e., we must have $v_0=0$ and $\v \perp \hat \p$. This is exactly the 
definition of $\Pi_{\p,\s}$ in the case $\hat\s = -\hat \p$ under consideration.
\qed

\begin{rmk}
When $\hat\s\neq-\hat\p$, the circle $\mathscr{C}_{\p,\s}$ can be explicitly 
described in terms of $\tau = \sqrt{|\s|/|\p|}$, $\u_{\p,\s}$, and a parameter 
$\phi$ as
\begin{equation}
\mathscr{C}_{\p,\s} = 
\left\{\, \tau \,\u_{\p,\s}\, (\cos\phi+\sin\phi\,\hat\p) \,\left|\, 
\phi\in[\,0,2\pi) \right. \,\right\}.
\end{equation}
\end{rmk}

Using again the notations $\X=x_0+\X$, $\P=p_0+\p$, $\Q=q_0+\q$, 
$\R=r_0+\r$, $\S=s_0+\s$ for the splittings of the variable and coefficients 
into scalar and vector parts, we are now ready to complete the study of 
the general equation.

\begin{thm}
\label{thm:nonreal}
Defining $\tilde\R$ and $\tilde\S=(\tilde{s}_0,\tilde\s)$ by
\begin{align*}
\tilde \R &:= (\Q^*\P-\R\P^*) (\P-\P^*)^{-1} \,, \\ 
\tilde \S &:= \S + \frac{\Q^*\P\R^*+ \R\P^*\Q}{|\P-\P^*|^{2}} + 
\frac{\R(\P^*-\P)\Q}{|\P-\P^*|^{2}} - \frac{\Q^*\P^*\Q}{|\P-\P^*|^{2}} 
- \frac{\R\P^*\R^*}{|\P-\P^*|^{2}} \,,
\end{align*}
the solutions in $\hh$ to equation~\eqref{Xeqn} with $\P,\Q,\R,\S \in \hh$ and 
$\P \not\in \rr$ are given by
\[
\X = \Z - (\R-\Q^*)(\P-\P^*)^{-1} ,
\]
where $\Z \in \mathscr{C}_{\p,\tilde \s}$ satisfies
\be
\label{hyperplane}
\langle\Z,\tilde\R\rangle \,=\, \frac{\tilde s_0|\p|-p_0|\tilde\s|}{2|\p|}.
\ee
Specifically, the following cases may be distinguished.
\begin{enumerate}
\item 
If $\tilde\S=0$, then $\X=-(\R-\Q^*)(\P-\P^*)^{-1}$ is the unique solution.
\item 
If $\tilde\S$ is a non--zero real number, there is no solution. 
\item 
If $\tilde \S$ is not a real number, the following cases arise. 
\begin{enumerate}
\item 
If $\tilde\R \perp \Pi_{\p,\tilde \s}$ the set of solutions $\X$ is either 
the circle 
\[
\mathscr{C}_{\p,\tilde \s}- (\R-\Q^*)(\P-\P^*)^{-1}
\]
or the empty set, depending on whether or not $\tilde s_0|\p|=p_0|\tilde \s|$.
\item 
If $\tilde\R \not\perp \Pi_{\p,\tilde \s}$ there are two, one or no 
solutions, namely, the points $\X=\Z-(\R-\Q^*)(\P-\P^*)^{-1}$ where 
$\Z$ is an intersection point of the circle $\mathscr{C}_{\p,\tilde \s}$ 
with the coplanar affine line specified by intersecting $\Pi_{\p,\tilde\s}$ 
with the hyperplane \eqref{hyperplane}. 
\end{enumerate}
\end{enumerate}
\end{thm}

\proof
Upon setting $\Z:=\X + (\R-\Q^*)(\P-\P^*)^{-1}$ and substituting $\X=\Z-
(\R-\Q^*)(\P-\P^*)^{-1}$ into equation \eqref{Xeqn}, we obtain
\begin{align*}
\S =\ & \Z\P\Z^*- \Z \P (\P^*-\P)^{-1}(\R^*-\Q) + \Z \Q \\
&- (\R-\Q^*)(\P-\P^*)^{-1} \P \Z^* + \R \Z^*\\
&+ (\R-\Q^*)(\P-\P^*)^{-1} \P (\P^*-\P)^{-1}(\R^*-\Q)\\
&- (\R-\Q^*)(\P-\P^*)^{-1} \Q - \R (\P^*-\P)^{-1}(\R^*-\Q)\\
=\ & \Z\P\Z^*+\Z  (\P^*-\P)^{-1}[-\P\R^*+ \P\Q + (\P^*-\P)\Q] \\
&+ [\Q^*\P + \R (-\P + \P-\P^*)] (\P-\P^*)^{-1}\Z^*\\
&+ |\P-\P^*|^{-2}(\R-\Q^*){\P}(\R^*-\Q)\\
&-  |\P-\P^*|^{-2}(\R-\Q^*)({\P^*-\P})\Q \\
&-  |\P-\P^*|^{-2}\R ({\P-\P^*})(\R^*-\Q)\\
=\ & \Z\P\Z^*+\Z  (\P^*-\P)^{-1}(-\P\R^*+\P^*\Q)\\
&+ (\Q^*\P-\R\P^*) (\P-\P^*)^{-1}\Z^*\\
&+ |\P-\P^*|^{-2}(\R-\Q^*)(\P\R^*-\P^*\Q) \\
&+  |\P-\P^*|^{-2} \R(-\P\R^*+\P\Q+\P^*\R^*-\P^*\Q)\\
=\ & \Z\P\Z^*+\Z \tilde \R^*+ \tilde \R \Z^*\\
&+ |\P-\P^*|^{-2}[-\Q^*(\P\R^*-\P^*\Q)+ \R(-2\P^*\Q+\P\Q + \P^*\R^*)] ,
\end{align*}
where $\tilde\R=(\Q^*\P-\R\P^*)(\P-\P^*)^{-1}$. This gives
\be\label{nonrealreduction}
\Z\P\Z^*+ \Z \tilde \R^* + \tilde \R \Z^* = \tilde \S
\ee
where
\begin{align*}
\tilde \S =\ & \S + \frac{\Q^*(\P\R^*-\P^*\Q)+ \R(2\P^*\Q 
-\P\Q -\P^*\R^*)}{|\P-\P^*|^{2}} \\
=\ & \S + \frac{\Q^*\P\R^*+ \R\P^*\Q}{|\P-\P^*|^{2}} + 
\frac{\R(\P^*-\P)\Q}{|\P-\P^*|^{2}} - \frac{\Q^*\P^*\Q}{|\P-\P^*|^{2}} 
- \frac{\R\P^*\R^*}{|\P-\P^*|^{2}}.
\end{align*}
Equation~\eqref{nonrealreduction} is equivalent to the system
\be
\label{nonrealsystem}
p_0\, |\Z|^2 + \Z \tilde \R^* + \tilde \R \Z^* = \tilde s_0 , \qquad
\Z\p\Z^* \,=\, \tilde \s .
\ee
If $\tilde\s=0$, the only solution to the second of these equations is 
$\Z=0$.  This also satisfies the first equation if and only if $\tilde s_0
=0$. We have thus established  cases {\it 1} and {\it 2}. For case {\it 3}, 
with $\tilde\s\neq0$, equations~\eqref{nonrealsystem} are equivalent to
\be
\label{nonrealsystem2}
\Z \tilde \R^* + \tilde \R \Z^* = 
\frac{\tilde s_0|\p|-p_0|\tilde \s|}{|\p|} , \qquad
\Z\p\Z^* \,=\, \tilde \s .
\ee
By Lemma~\ref{circle}, the solutions of the latter equation comprise the 
circle $\mathscr{C}_{\p,\tilde \s}$ of radius $\sqrt{|\tilde \s|/|\p|}$ and 
center $0$ in the $2$--plane $\Pi_{\p,\tilde\s}$.  The solutions to the 
system \eqref{nonrealsystem2} are determined by intersecting this circle 
with the set $\mathscr{H}$ of solutions to
\[
\langle \Z,\tilde\R \rangle \,=\, 
\frac{\tilde s_0|\p|-p_0|\tilde\s|}{2\,|\p|} .
\]
For any $\Z_0\in \mathscr{H}$, this equation is equivalent to 
\[
\langle \Z-\Z_0,\tilde\R \rangle \,=\, 0.
\]
Hence, the set $\mathscr{H}$ is an affine space orthogonal to $\tilde \R$. 
Case {\it 3(a)} is verified by observing that, since $\mathscr{C}_{\p,
\tilde \s}$ is a circle centered at $0$ in the plane $\Pi_{\p,\tilde\s}$, 
the following are equivalent:
\bi
\item 
$\mathscr{C}_{\p,\tilde \s}$ is the set of solutions to system
\eqref{nonrealsystem2};
\item 
$\Pi_{\p,\tilde \s}\subset \mathscr{H}$;
\item 
$0 \in \mathscr{H}$ and  the equality $\langle \Z,\tilde\R \rangle \,=\, 0$ 
holds for all $\Z\in \Pi_{\p,\tilde \s}$; 
\item 
$\tilde s_0\,|\p|=p_0\,|\tilde\s|$ and $\tilde\R \perp \Pi_{\p,\tilde \s}$.
\ei
Finally, consider case {\it 3(b)} with $\tilde\R \not\perp \Pi_{\p,\tilde \s}$ 
(whence $\tilde\R \neq 0$). Then $\mathscr{H}$ is a hyperplane that, by 
Lemma~\ref{circle}, intersects $\Pi_{\p,\tilde \s}\supset\mathscr{C}_
{\p,\tilde\s}$ in the affine line defined by the equations
\[
\Z\,\hat\p-\hat {\tilde {\s}}\, \Z = 0 , \qquad
\Z \tilde \R^* + \tilde \R \Z^* =  
\frac{\tilde s_0|\p|-p_0|\tilde \s|}{|\p|} .
\]
Finally, the translation $\X= \Z-(\R-\Q^*)(\P-\P^*)^{-1}$ gives the solutions
of quation~\eqref{Xeqn} in terms of those of~\eqref{nonrealreduction}.
\qed

\bigskip

Cases {\it 3(a)} and {\it 3(b)} of Theorem~\ref{thm:nonreal} are 
characterized by whether or not $\tilde\R$ is orthogonal to the plane 
$\Pi_{\p,\tilde\s}$. In an algorithm, this could be determined by checking 
to see if the scalar product of $\tilde\R$ with any two linearly independent 
quaternions in $\Pi_{\p,\tilde\s}$ vanishes. Note also that the computation 
of the points $\Z$ in case {\it 3(b)} can be performed by: (a) writing a 
parameterization $\xi \mapsto \A\, \xi + \B$
of the affine line specified by intersecting $\Pi_{\p,\tilde\s}$ 
with the hyperplane \eqref{hyperplane}; and (b) solving the real quadratic equation $| \A\, \xi + \B|^2=|{\tilde {\s}}|^2/|{\tilde {\p}}|^2$.

%%%%%%%%%%%%%%%%%%%%%%%%%%%%%%%%%%%%%%%%%%%%%%%%%%%%%%%%%%%%%%%%%%%%%%%%%%%%

\medskip\noindent

We now illustrate cases {\it 1}. and {\it 2}. of Theorem~\ref{thm:nonreal} 
by the following examples.

\begin{exm}
{\rm If $\P = \i$, $\Q=\R=1$, and $\S=0$, equation~\eqref{Xeqn} becomes
\[
\X\,\i\,\X^*+\X+\X^*=0,
\]
which has the unique solution $\X=0$ in $\hh$.}
\end{exm}

\begin{exm}
{\rm If $\P = \i$ and $\Q=\R=\S=1$, equation~\eqref{Xeqn} becomes
\[
\X\,\i\,\X^*+\X+\X^*=1 ,
\]
which has no solution in $\hh$.}
\end{exm}

\medskip\noindent
The following two instances exemplify case {\it 3(a)} in Theorem~\ref
{thm:nonreal}.

\begin{exm}
{\rm If $\P= \i$, $\Q=\R=1$, and $\S=-\i$, equation~\eqref{Xeqn} becomes
\[
\X\,\i\,\X^*+\X+\X^*= -\i.
\]
in this case, the set of solutions in $\hh$ is the circle
\[
\mathscr{C}_{\i,-\i} = \{\, x_2 \j+ x_3 \k \,|\, x_2^2+x_3^2=1 \,\}.
\]
}
\end{exm}

\begin{exm}
{\rm If $\P= \i$, $\Q=\R=1$, and $\S=1-\i$, equation~\eqref{Xeqn} becomes
\[
\X\,\i\,\X^*+\X+\X^*= 1-\i,
\]
which is impossible to satisfy, since the (purely imaginary) circle 
$\mathscr{C}_{\i,-\i}$ does not intersect the hyperplane $\X+\X^*=1$.}
\end{exm}

\medskip\noindent
Case {\it 3(b)} of Theorem~\ref{thm:nonreal} is illustrated by the 
following family of examples.

\begin{exm}
{\rm If $\P=\i$, $\Q=\R=1$, and $\S=s_0+\i$, equation~\eqref{Xeqn} becomes
\[
\X\,\i\,\X^*+\X+\X^*= s_0+\i.
\]
The solutions are the intersections of the circle $\mathscr{C}_{\i,\i} 
=\{\cos\phi+\sin\phi\,\i : \phi \in \rr\}$ with the hyperplane $\X+\X^*
=s_0$. If $s_0 = 2\cos\phi_0$ there are two or one solutions in $\hh$, 
namely $\cos\phi_0\pm\sin\phi_0\,\i$. If $|s_0|>2$, on the other hand, 
there is no solution in $\hh$.}
\end{exm}

Further examples of numerical solutions to equation~\eqref{Xeqn}, 
computed in the context of the surface construction problem (Section~\ref
{sec:phcurves}), are presented in \cite{farouki15}.

%%%%%%%%%%%%%%%%%%%%%%%%%%%%%%%%%%%%%%%%%%%%%%%%%%%%%%%%%%%%%%%%%%%%%%%%%%%%

\section{Closure}
\label{sec:close}

Although the solution of equations in the space of quaternions $\hh$ 
has recently attracted considerable attention, most studies have been 
restricted to the case of unilateral coefficients. In the present study, we 
have considered a special quadratic quaternion equation in the quaternion 
variable \emph{and} its conjugate, with mixed coefficients. This equation,
arising from a surface construction problem \cite{farouki15}, was shown to 
admit a complete characterization of its solutions, for all possible instances 
of the coefficients. In addition to point solutions, \emph{circles} or 
\emph{3--spheres} of solutions are observed --- as distinct from the case of 
unilateral coefficients, which admits \cite{gentili08b,pogorui04,topuridze09} 
only point solutions and 2--spheres of solutions. We have thus determined 
a significant class of low--degree quaternion equations (disjoint from 
that of polynomial equations with unilateral coefficients), for which a 
comprehensive solution can be achieved.

%%%%%%%%%%%%%%%%%%%%%%%%%%%%%%%%%%%%%%%%%%%%%%%%%%%%%%%%%%%%%%%%%%%%%%%%%%%%

\subsection*{Acknowledgements}
{\small
This work was supported by the following grants of the Italian Ministry of Education (MIUR): Futuro in Ricerca {\it DREAMS} (RBFR13FBI3); Futuro in Ricerca {\it Differential Geometry and Geometric Function Theory} (RBFR12W1AQ); and PRIN {\it Variet\`a reali e complesse: geometria, topologia e analisi armonica} (2010NNBZ78).
It was also supported by the following research groups of the Istituto Nazionale di Alta Matematica (INdAM): Gruppo Nazionale per il Calcolo Scientifico (GNCS) and Gruppo Nazionale per le Strutture Algebriche, Geometriche e le loro Applicazioni (GNSAGA).
}
%%%%%%%%%%%%%%%%%%%%%%%%%%%%%%%%%%%%%%%%%%%%%%%%%%%%%%%%%%%%%%%%%%%%%%%%%%%%

\def\AML{{\it Appl.\ Math.\ Lett.\ }}
\def\AMM{{\it Amer.\ Math.\ Monthly\ }}
\def\ACM{{\it Adv.\ Comp.\ Math.\ }}
\def\ACMTMS{{\it ACM Trans.\ Math.\ Software\ }}
\def\ACMTOG{{\it ACM Trans.\ Graphics\ }}
\def\BAMS{{\it Bull.\ Amer.\ Math.\ Soc.\ }}
\def\CA{{\it Comm.\ Alg.\ }}
\def\CAD{{\it Comput.\ Aided Design\ }}
\def\CAEJ{{\it Comput.\ Aided Eng. J.\ }}
\def\CAGD{{\it Comput.\ Aided Geom.\ Design }}
\def\CAVW{{\it Comput.\ Anim.\ Virt.\ Worlds }}
\def\CG{{\it Computers \& Graphics }}
\def\CGIP{{\it Comput.\ Graphics Image\ Proc.\ }}
\def\CMA{{\it Comput.\ Math.\ Applic.\ }}
\def\CV{{\it Complex Var.\ }}
\def\CVGIP{{\it Comput.\ Vision, Graphics, Image\ Proc.\ }}
\def\EJLA{{\it Electron.\ J.\ Lin.\ Alg.\ }}
\def\GM{{\it Graph.\ Models\ }}
\def\IBMJRD{{\it IBM J.\ Res.\ Develop.\ }}
\def\JCAM{{\it J.\ Comput.\ Appl.\ Math.\ }}
\def\IEEECGA{{\it IEEE Comput. Graph. Applic.\ }}
\def\IJAMT{{\it Int.\ J.\ Adv.\ Manuf.\ Tech.\ }}
\def\IJMS{{\it Int.\ J.\ Model.\ Sim.\ }}
\def\IJMTM{{\it Int.\ J.\ Mach.\ Tools Manuf.\ }}
\def\IJRR{{\it Int.\ J.\ Robot.\ Res.\ }}
\def\IMAJNA{{\it IMA J.\ Numer.\ Anal.\ }}
\def\IUJM{{\it Ind.\ Univ.\ J.\ Math.\ }}
\def\JMS{{\it J.\ Math.\ Sci.\ }}
\def\JC{{\it J.\ Complexity\ }}
\def\JSC{{\it J.\ Symb.\ Comput.\ }}
\def\MC{{\it Math.\ Comp.\ }}
\def\MMA{{\it Math.\ Model.\ Anal.\ }}
\def\MMAS{{\it Math.\ Methods\ Appl.\ Sci.\ }}
\def\MJM{{\it Milan J.\ Math.\ }}
\def\MMJ{{\it Mich.\ Math.\ J.\ }}
\def\MZ{{\it Math.\ Z.\ }}
\def\NA{{\it Numer.\ Algor.\ }}
\def\NMTMA{{\it Numer.\ Math.\ Theor.\ Meth.\ Appl.\ }}
\def\PAMS{{\it Proc.\ Amer.\ Math.\ Soc.\ }}
\def\SIAMJNA{{\it SIAM J.\ Numer.\ Anal.\ }}
\def\SIAMR{{\it SIAM Rev.\ }}
\def\TAMS{{\it Trans.\ Amer.\ Math.\ Soc.\ }}

\begin{flushleft}

\end{flushleft}

\end{document}